\documentclass[letterpaper, 11 pt]{amsart}
\usepackage[short]{srollens-en}
\usepackage{color}

\title[Torsion and cotorsion in the sheaf of K\"ahler  differentials]{Torsion and cotorsion in the sheaf of K\"ahler  differentials on some mild singularities}
\date{\today}
\author{Daniel Greb}
\address{Daniel Greb\\Institut f\"ur Mathematik\\Abteilung f\"ur Reine Mathematik\\Albert--Ludwigs--Universit\"at Freiburg\\Eckerstra{\ss}e 1\\79104 Freiburg im Breisgau\\Germany}
\curraddr{Mathematics Department\\ Princeton University\\ Fine Hall, Washington Road\\
Princeton NJ 08544-1000, USA}

\email{daniel.greb@math.uni-freiburg.de}

\author{S\"onke Rollenske}
\address{S\"onke Rollenske\\
Institut f\"ur Mathematik\\
Johannes Gutenberg Universit\"at Mainz\\
55099 Mainz\\
Germany}
\email{rollensk@uni-mainz.de}

\DeclareMathOperator{\Tors}{\mathrm{Tors}}
\DeclareMathOperator{\Sym}{\mathrm{Sym}}

{\theoremstyle{Lehn-up}
\newtheorem{que}[theo]{Question}
}

\newcommand{\DDual}[1]{{#1}^{\vee\vee}}

\begin{document}

\begin{abstract}
We give a criterion for the sheaf of K\"ahler differentials on a cone over a smooth projective variety to be torsion\-free.

Applying this to Veronese embeddings of projective space and using known results about differentials on quotient singularities we show that even for mild,~e.\,g.~Gorenstein  terminal, singularities the sheaf of K\"ahler differentials will in general have torsion and cotorsion.
\end{abstract}


\maketitle

\section{Introduction}
Let $Z$ be an algebraic variety over a field $k$, which  we assume to be of characteristic 0.
One of the few objects that come naturally with $Z$ is its sheaf of K\"ahler differentials $\Omega_Z=\Omega_{Z/k}$. It is the sheafified version of the module of K\"ahler differentials, which is also an important tool in commutative algebra. In this note we stick to the geometric language.

The sheaf of differentials and its higher exterior powers play an important r\^ole in  many contexts, most prominently deformation theory, vanishing theorems {and} (if $Z$ is sufficiently nice) duality theory. For a more local example one could mention Berger's conjecture that a curve is smooth if and only if its sheaf of differentials is torsion\-free (see e.g. \cite{berger63, GreuelDeligne, pohl91}) or the Zariski--Lipman conjecture: if $\Omega_Z^\vee$ is locally free then $Z$ is regular (see e.g. \cite{platte88, blls02}).

In the context of the minimal model program Greb, Kebekus, Kov\'acs and Peternell  proved strong extension theorems for differential forms \cite{gkkp10}, but instead of the sheaf of K\"ahler differentials itself they used  its reflexive hull $\DDual{\Omega_Z}$ also called module of Zariski differentials (see e.g. \cite{knighten73}). For applications the following obvious question comes to mind:
\begin{que}\label{quest}
 If $Z$ has \emph{mild} singularities, is $\Omega_Z$ reflexive or at least torsion\-free?
\end{que}
Phrasing this slightly differently we ask if the natural map $\phi\colon\Omega_Z\to \DDual{\Omega_Z}$, whose kernel is the torsion submodule  $\Tors(\Omega_Z)$,  is bijective or at least injective. In the terminology of \cite[(1.6)]{ReidYPG} we say that $\Omega_Z$ has \emph{cotorsion} if $\phi$ is not surjective.\footnote{{In the case of curve singularities torsion and cotorsion in $\Omega$ were extensively studied for example by Greuel and his collaborators in
\cite{GreuelBuchweitz} and \cite{GreuelPfister}.}}

If one translates \emph{mild singularity} as being a local complete intersection, then indeed a complete answer to Question~\ref{quest} is known.
\begin{theo}[\cite{Kunz}, Proposition 9.7, Corollary 9.8] \label{lci}
Let $Z$ be a local complete intersection. Then $\Omega_Z$ satisfies Serre's condition $S_d$ if and only if $Z$ is regular in codimension $d$.

In particular, if $Z$ is a normal local complete intersection, then $\Omega_Z$ is torsion\-free, and it is reflexive if and only if $Z$ is is nonsingular in codimension 2.
\end{theo}

However, in the context of modern birational geometry one usually measures the singularities of a normal variety in terms of discrepancies, which give rise to the definition of terminal, canonical, log terminal and other singularities \cite[Section 2.3]{Kollar-Mori}. Even terminal singularities, the mildest class considered, are in general not complete intersections, so {Theorem~\ref{lci}} does not apply.

Somewhat contrary to our expectations we will show below that the answer to Question \ref{quest} is essentially negative if one interprets mild in the sense of birational geometry; as soon as one leaves the world of local complete intersections one should expect the sheaf of K\"ahler differentials to have both torsion and cotorsion.

There are two cheap ways to produce non--lci singularities: quotients of finite groups and affine cones over projective varieties. In the second case the algebraic description is somewhat simpler and we give a criterion for the existence of torsion differential in Section \ref{criterion}. Recall that if $X\subset \IP^n$ is a projective variety or scheme  given by an  ideal sheaf $\ki$ then $X$ is called \emph{projectively normal} if $H^0(\IP^n, \ko_{\IP^n}(d))\onto H^0(X, \ko_X(d))$ for all $d\geq0$ or, equivalently, $H^1(\IP^n, \ki(d))=0$ for all $d\geq 0$.

\begin{theo}\label{main}
Let $X\subset \IP^n_k$ be a smooth projective variety over a field $k$ of characteristic zero and let $\ki$ be the sheaf of ideals defining $X$. Let $C_X\subset \IA^{n+1}$ be the affine cone over $X$.

If $H^1(\IP^n, \ki^2(d))=0$ for all $d\geq 0$ then $\Omega_{C_X}$ is torsion\-free.

If in addition $X$ is projectively normal, then $\Omega_{C_X}$ is torsion\-free if and only if also  the first infinitesimal neighbourhood of $X$ {in $\IP^n_k$} is projectively normal.\footnote{{After the publication of the first version of this preprint J. Wahl brought to our attention that he had already observed this in the projectively normal case \cite[Proposition 1.4]{wahl97}.}}
\end{theo}

In the cone situation the study of  cotorsion is more problematic but in Section \ref{cotorsion} we recall results by Knighten and Steenbrink that give an easy sufficient criterion for cotorsion on finite quotient singularities.

In the last section we study our main class of examples, namely $X_{r,d}$, the affine cone over the $d$th Veronese embedding of $\IP^r$. We can also describe $X_{r,d}$ as a cyclic quotient singularity.  Such cones have torsion differentials if and only if $d\geq3$ (Proposition \ref{Vrd}) and cotorsion as soon as $d\geq 2$ (Proposition \ref{Vrdkt}); we have collected some significant cases in Table \ref{tab}\footnote{The sheaf of K\"ahler differentials for $X_{1,2}$ and $X_{1,3}$ was also computed in \cite{knighten73} but the examples were not widely known.}.
\begin{table}[ht]\label{tab}
\caption{Some Veronese cones with torsion or cotorsion in $\Omega_{X_{d,r}}$.}
\begin{center}
\begin{tabular}{ccccccc}
singularity & $\dim$ & type & Gorenstein & torsion & cotorsion \\
 \hline
$X_{1,2}$ (A1)& 2  & canonical & yes & no & yes\\
$X_{1,3}$ & 2  & log terminal & no & yes & yes\\
$X_{2,2}$ & 3 & terminal & no & no & yes\\
$X_{2,3}$ & 3 & canonical & yes & yes & yes\\
$X_{3,2}$ & 4 & terminal & yes & no & yes\\
$X_{3,3}$ & 4 & terminal & no & yes & yes\\
$X_{5,3}$ & 6 & terminal & yes & yes & yes\\
\end{tabular}
\end{center}
\end{table}

\begin{rem}
In the surface case our results are optimal in the following sense: let $Z$ be a surface singularity. If $Z$ is terminal then it is smooth and thus $\Omega_Z$ is locally free. If $Z$ is canonical but not terminal then it is one of the well--known ADE singularities, thus a hypersurface singularity; by Theorem~\ref{lci} the sheaf of K\"ahler differentials $\Omega_Z$ is torsion\-free but not reflexive in this case. The easiest log terminal point is the cone over the twisted cubic $X_{1,3}$ and in this case $\Omega_Z$ has both torsion and cotorsion.

{Gorenstein terminal 3--fold singularities are hypersurface singularities by a result of Reid \cite[(3.2) Theorem]{ReidYPG}; hence Theorem~\ref{lci} applies to show that the sheaf of K\"ahler differentials is reflexive in this case.}

It is possible that other classes of mild singularities in small dimensions turn out to have torsion\-free or reflexive sheaf of K\"ahler differentials as well. Using structural results like Hilbert--Burch or Buchsbaum--Eisenbud, low--codimensional singularities might also be accessible (compare \cite[Section 4]{mvs01}).
\end{rem}

\section{Torsion differentials on cones}\label{criterion}
In this section we give the proof of Theorem~\ref{main}. Let $X\subset \IP^n$ be a smooth, irreducible, non--degenerate projective variety and let $C_X$ be the affine cone over $X$. We denote by $\ki$ the ideal sheaf of $X$.

The question {whether} $\Omega_{C_X}$ has torsion is purely algebraic, because the cone is affine. Denote by $S=k[x_0, \dots, x_n]$ the polynomial ring with homogeneous maximal ideal $\gothm=(x_0, \dots, x_n)$ and by $R=S/I$ the homogeneous coordinate ring of $X$ whose homogeneous maximal ideal we denote by $\gothn$.

Note that $C_X$ is smooth outside the vertex so we are only interested in the local behaviour at the vertex. In the following we will use some facts about local cohomology all of which can be found in \cite[Appendix~1]{EisenbudSyz} or in {more} detail in \cite{HartshorneLC}.

In the following diagram the middle column is the conormal sequence and the maps to the maximal ideal are given by the contraction with the vectorfield $\xi=\sum_ix_i\frac \del{\del x_i}$.
\[\xymatrix{
& I/I^2\ar[d] \ar@{=}[r] & I/I^2\ar[d]\\
0\ar[r]& M \ar[r]\ar[d] &\Omega_S\tensor R \ar[d]\ar[r]^{\xi\lrcorner}& \gothn \ar[r]\ar@{=}[d] & 0\\
0 \ar[r]& N\ar[r]\ar[d] & \Omega_R \ar[r]^{\xi\lrcorner}\ar[d] &\gothn \ar[r] &0\\
& 0& &}
\]
{By definition,} the modules $M$ and $N$ are the respective kernels {of the contraction maps $\xi\lrcorner$.} The fact that the composition $I/I^2\to \Omega_S\tensor R\to \gothn$ is zero follows from Euler's formula: for every homogeneous element $f\in I$ we have $f=\sum_i x_i\frac{\del f}{\del x_i}$.

Since $\gothn$ is torsion\-free the torsion--submodule $\Tors(\Omega_R)$ {is isomorphic to} $\Tors(N)$. Since $C_X$ is smooth outside the vertex we can compute {this} torsion submodule via local cohomology
\[ \Tors(N)=H^0_\gothn(N),\]
which we now relate to sheaf cohomology on $\IP^n$.

Introducing the appropriate grading we can transform the above diagram into a diagram of coherent sheaves on $\IP^n$, where the middle row is the restriction of the Euler Sequence to $X$ and the first column becomes the conormal sheaf sequence for $X$. {The map $\ki/\ki^2 \to \Omega_{\IP^n}\restr X $ in the latter sequence} is injective because $X$ is smooth.
\[\xymatrix{
& 0\ar[d]& 0 \ar[d]\\
& \ki/\ki^2\ar[d] \ar@{=}[r] & \ki/\ki^2\ar[d]\\
0\ar[r]& \Omega_{\IP^n}\restr X \ar[r]\ar[d] &\ko_X(-1)^{\oplus n+1} \ar[d]\ar[r]& \ko_X\ar[r]\ar@{=}[d] & 0\\
0 \ar[r]& \Omega_X\ar[r]\ar[d] & \widetilde{\Omega_R} \ar[r]\ar[d] &\ko_X \ar[r] &0\\
& 0& 0
}
\]

We now use the comparison sequence for local and sheaf cohomology (\cite[Cor.~A1.12]{EisenbudSyz}) for the first column of the diagram {to obtain}
\begin{equation}\label{diag}
\xymatrix
{ &&& 0 \ar[d]\\
&0\ar[r] & I/I^2 \ar[r] \ar[d] & \Gamma_*(\ki/\ki^2) \ar[r] \ar[d] & H^1_\gothn(I/I^2)\ar[r]\ar[d] &0\\
&0\ar[r]& M \ar[r]\ar[d] &\Gamma_*(\Omega_{\IP^n}\restr X)\ar[d]\ar[r] & H^1_\gothn(M)\ar[r]\ar[d] &0\\
0\ar[r] & H^0_\gothn(N)\ar[r] & N \ar[r] \ar[d] & \Gamma_*(\Omega_X) \ar[r] \ar[d] & H^1_\gothn(N)\ar[r]&0\\
 &&0&
}
\end{equation}
where we already used the following Lemma.
\begin{lem} In the above situation we have
 \begin{enumerate}
  \item $H^0_\gothn(I/I^2)=0$,
\item $H^0_\gothn(M)=0$.
 \end{enumerate}
\end{lem}
\begin{proof}
\begin{enumerate}
\item This follows from the short exact sequence $0\to I^2\to I \to I/I^2\to 0$ and the fact that $I$ and $I^2$ are saturated and have depth at least 2.
\item $M$ is a submodule of the free module $\Omega_S\tensor R$, hence torsion\-free.
\end{enumerate}
\end{proof}

Applying the same arguments as in the proof of the snake lemma to \eqref{diag} we get an exact sequence
\begin{equation}\label{snake} 0\to  H^0_\gothn(N)\to  H^1_\gothn(I/I^2)\to H^1_\gothn (M) \to H^1_\gothn (N)  \end{equation}
and thus have shown that $\Omega_X$ is torsion\-free if $ H^1_\gothn(I/I^2)=0$.
This latter group can be interpreted as follows: the local cohomology sequence for $0\to I^2\to I \to I/I^2\to 0$ compared with the sheaf cohomology sequence of $0\to \ki^2\to \ki\to \ki/\ki^2\to 0$ yields a diagram
\begin{equation}\label{iki}
\xymatrix{
 I \ar[r]\ar[d]^\isom & I/I^2 \ar[r]\ar@{^(->}[d] &0 \ar[d]\\
\Gamma_*(\ki) \ar[r]& \Gamma_*(\ki/\ki^2)\ar[r]\ar@{->>}[d] & \bigoplus_d H^1(\ki^2(d))\ar[r]&\bigoplus_d H^1(\ki(d))\\
&  H^1_\gothn(I/I^2)\ar@{^(->}[ru]
}
\end{equation}
We have thus proved the first part of  Theorem \ref{main} from the introduction which we repeat here for the convenience of the reader.
\begin{custom}[Theorem \ref{main}]
Let $X\subset \IP^n_k$ be a smooth projective variety over a field $k$ of characteristic 0 and let $\ki$ be the sheaf of ideals defining $X$. Let $C_X\subset \IA^{n+1}$ be the affine cone over $X$.

If $H^1(\IP^n, \ki^2_X(d))=0$ for all $d\geq0$ then $\Omega_{C_X}$ is torsion\-free.

If in addition $X$ is projectively normal, then $\Omega_{C_X}$ is torsion\-free if and only if also  the first infinitesimal neighbourhood of $X$ is projectively normal.
\end{custom}
\begin{proof}[Proof of the second part]
By definition, projective normality {of $X$} is equivalent to $H_\gothn^0(R)=H_\gothn^1(R)=0$ or equivalently $\bigoplus_d H^1(\IP^n, \ki(d))=0$.
The exact sequence in local cohomology
\[0=H^0_\gothn(\gothn) \to H^1_\gothn(M)\to H^1_\gothn(\Omega_S\tensor R)=H^1_\gothn(R)^{\oplus n+1}=0, \]
induced by $0\to M\to \Omega_S\tensor R\to \gothn\to 0$, shows that $H^1_\gothn(M)=0$.
Together with \eqref{snake} and \eqref{iki} {this vanishing implies} that the composition
\[ \Tors(\Omega_R)=H^0_\gothn(N)\into H^1_\gothn(I/I^2)\into \bigoplus_dH^1(\IP^n, \ki^2(d))\]
is an isomorphism. {This establishes} the desired equivalence.
\end{proof}

\section{Zariski differentials on quotient singularities and cotorsion}\label{cotorsion}
We briefly recall a result of Knighten, also discovered by Steenbrink \cite[Lem.~1.8]{SteenbrinkHodgeVanishingCohomology}, {describing the double dual of $\Omega_{X/G}$ on finite quotient singularities $X/G$}, and deduce a sufficient condition for the existence of cotorsion in the sheaf of K\"ahler differentials.

Let $A$ be a regular local $k$--algebra and let $G$ be a finite group acting on $A$.
We denote by
$\pi\colon X=\spec A \to X/G=\spec A^G$
the quotient map.
\begin{theo}[\cite{knighten73}, Theorem 3]
The natural map $\DDual{\Omega_{X/G}}\to (\pi_*\Omega_X)^G$ is an isomorphism.
\end{theo}

In concrete situations the module of invariant differentials is not hard to compute and and the following will turn out to be useful
\begin{cor}\label{cotors}
Assume that $X/G$ has an isolated singularity.
Let $\gothm$ be the maximal ideal of $A^G$ and $e=\dim_k\gothm/\gothm^2$ be the embedding dimension of $X/G$. If the minimal number of generators for $(\Omega_X)^G$ (considered as a $A^G$--module) is bigger than $e$ then $\Omega_{X/G}$ has cotorsion and is not reflexive.
\end{cor}
\begin{proof}
By the conormal sequence the sheaf $\Omega_{X/G}$ of K\"ahler differentials can be generated by $e$ elements and thus can never surject onto $(\pi_*\Omega_X)^G$ if this {sheaf needs} more than $e$ generators.
\end{proof}

\section{Examples: Cones over Veronese embeddings}\label{examples}

In this section we study cones over Veronese embeddings, including the examples mentioned in Table \ref{tab}. The first subsection collects some general properties {of these cones}, torsion differentials {are computed} in Propostion \ref{Vrd}, and cotorsion will be discussed in Proposition \ref{Vrdkt}.

\subsection{Basic properties}\label{subsect:basicprop}
Here we collect basic properties of the cones over the Veronese embeddings. In particular, we discuss their realisation as cyclic quotient singularities, and compute the discrepancies of the canonical resolution in order to determine in which cases {these cones}  are terminal, canonical, etc. For this we follow \cite[\S 1]{reid80} and \cite[Sect.~7.2]{Debarre}, where one can also find the definition of {the singularities appearing in the Minimal Model Program.}

Let $\mu_d$ be the cyclic group of order $d$ and let $\rho\colon \mu_d\to GL_{r+1}(k)$ be the representation given by choosing a primitive  $d$th root of unity $\xi$ and sending a generator of $\mu_d$ to
\[\begin{pmatrix}\xi &   &  \\
                               & \ddots &  \\
                               &        & \xi \end{pmatrix}  \in GL_{r+1}(k).\]
Let $X_{r,d} = \IA_k^{r+1} / \mu_d$ be the resulting quotient singularity (in Reid's \cite[(4.2)]{ReidYPG} notation, this is a singularity of type $\frac{1}{d}(1, 1, \dots, 1)$). Since the ring $\IC[x_0, \dots, x_r]^{\mu_d}$ is generated by all monomials of total degree $d$ in the coordinates $x_0, \dots, x_n$, the quotient $X_{r,d}$ is isomorphic to the cone $C_{r,d}$ over the image $V_{r,d}$ of the $d$th Veronese embedding $v_d\colon\IP^r\to \IP^n$ (so $n=\binom{r+d}{d}-1$). This isomorphism is induced by the map
\[(x_0, \dots, x_n) \mapsto (\dots, \prod_{a_0 + \cdots + a_n = d} x_0^{a_0} \cdots x_n^{a_n}, \dots). \]
The blow--up $\pi\colon Y_{r,d} \to X_{r,d}$ of the origin in $X_{r,d}$ is smooth, isomorphic to the total space of the line bundle  $\mathcal{O}_{X_{r,d}}(-1)$.

One can check (see \cite[p.278]{reid80}) that the index of $X_{r,d}$ is the denominator of $(r+1)/d$, i.e., $\mathrm{index}{(X)}= \frac{d}{\text{g.c.d}(r+1,d)}$. In particular, the canonical divisor $K_{X_{r,d}}$ is $\IQ$--Cartier, and $X_{r,d}$ is Gorenstein if and only if $d$ divides $r+1$.
We next compute the discrepancy of the unique exceptional divisor $E\isom V_{r,d}$; for simplicity we suppress the indices in the notation. Since the canoncial divisor is $\IQ$--Cartier we can write
\[K_Y \sim \pi^*K_X + aE\]
for a rational number $a$.
The normal bundle of $E$ in $Y$ is $\mathcal{O}_V(-1) \cong \mathcal{O}_{\IP^{r}}(-d)$. By adjunction, the canonical divisor of $E$ is $\mathcal{O}_E(K_Y + E)$. Hence, by restricting to $E$ we obtain
\[-r-1 = (a+1)(-d), \text{ hence } a = \frac{r+1}{d} -1. \]
It follows that all Veronese cones $X_{r,d}$ are log terminal, and in addition that
$X_{r,d}$ is terminal (canonical) if and only if $r+1>d$ ($r+1 \geq d$).

\subsection{Torsion differentials}

We now apply Theorem \ref{main} to the Veronese cones.

\begin{prop}\label{Vrd}
 Let $X_{r,d}$ be the affine cone over $V_{r,d}\subset \IP^n$, the image of the $d$th Veronese embedding of $\IP^r$ (so $n=\binom{r+d}{d}-1$). Then $\Omega_{X_{r,d}}$ has torsion if and only if $d\geq 3$.
\end{prop}
Note that by the discussion in the previous section all cones $X_{r,d}$ have a description as group quotients. So torsion differentials occur in abundance also on quotient singularities.
\begin{proof}
We now fix $r$ and $d$ and denote by $\ki$ the ideal {sheaf} of the image of the Veronese embedding $V=v_d(\IP^r)$.
To avoid confusion we denote by $H$ a hyperplane on $\IP^n$, so that $\ko_{\IP^n}(mH)\restr V=\ko_{\IP^r}(md)$.

Since $V$ is projectively normal the second part of Theorem \ref{main} applies and we only have to check if there is an $m$ such that $H^1(\IP^n,\ki^2(mH))\neq 0$. By a result of Wahl \cite[Theorem 2.1]{wahl97} we have
\[ H^1(\ki^2_V(mH))=0 \qquad \text{for $m\neq 2$}\]
and thus we only need to consider the case $m=2$.

We will start by showing that $\Omega_{X_{r,d}}$ is torsion\-free if $d=2$. Recall that there is a Gaussian map
\begin{equation}\label{gauss}
\Phi_{\ko_{\IP^r}(d)}: \Lambda^2 H^0(\IP^r,\ko_{\IP^r}(d)) \to H^0({\IP^r}, \Omega^1_{\IP^r}(2d))
\end{equation}
symbolically given by $s\wedge t\mapsto sdt-tds$, which in our case is a homomorphism of $\mathrm{SL}(r+1)$--representations. Wahl showed that  $H^1(\ki^2_V(2H))=\ker \Phi$ \cite[Proposition 1.8]{wahl97}. 

Note that the representation on $H^0(\IP^r, \Omega^1_{\IP^r}(4))$ is irreducible, cf.~\cite[Section 2]{wahl97}. It follows by a straightforward computation that the map $\Phi_{\ko_{\IP^r}(2)}$ is non-trivial and hence surjective. A calculation similar to the ones in Lemma \ref{lem:estimates} below gives the equality $\dim \Lambda^2 H^0(\IP^r,\ko_{\IP^r}(2))= h^0({\IP^r}, \Omega^1_{\IP^r}(4))$. Consequently, $\Phi_{\ko_{\IP^r}(2)}$ is an isomorphism. Thus $\Tors(\Omega_{X_{r,2}})=\ker(\Phi_{\ko_{\IP^r}(2)})=\{0\}$.
 This proves the ''only if'' part of the claim.

For the existence of torsion--differentials in case $d\geq3$ we give an elementary dimension estimate that does not depend on Wahls results.

By projective normality we have an exact sequence
\begin{multline}\label{seq1}
0\to H^0(\IP^n, \ki^2(2H))\to H^0(\IP^n, \ki(2H))\to \\
\to H^0(\IP^n, \ki/\ki^2(2H))\to H^1(\IP^n, \ki^2(2H))\to 0,
\end{multline}
and $H^1(\IP^n,\ki^2(2H))$ does not vanish if  $h^0(\IP^n, \ki/\ki^2(2H))>h^0(\IP^n, \ki(2H))$. This is the content of Lemma~\ref{lem:estimates} below.
\end{proof}
\begin{lem}\label{lem:estimates} In the situation above the following holds
\begin{enumerate}
 \item $h^0(\IP^n,\ki(mH))=\binom{\binom{d+r}r+m-1}{m} - \binom{md+r}r$.
\item $h^0(\IP^n, \ki/\ki^2(mH)) \geq \binom{d+r}r\binom{(m-1)d+r}r - \binom{md+r}r-(dm-1)\binom{r+dm-1}{md}$
\item $h^0(\IP^n, \ki/\ki^2(2H))-h^0(\IP^n, \ki(2H))>0$ for all $r\geq 1$, $d\geq 3$.
\end{enumerate}
\end{lem}
\begin{proof} If we dentote the graded polynomial ring in $r+1$ variables as $S=\bigoplus_{d\geq 0} S_d$ then the Veronese embedding is induced by the homomorphism of graded rings $\Sym^*(S_d) \to S$ with kernel a graded ideal $I$. In degree $m$ we get
\[0 \to H^0(\ki(mH))=I_m \to \Sym^m(S_d)\to S_{md}\to 0,\]
where surjectivity of the last map can easily be checked on monomials. The formula for the dimension follows from the well known formula for the dimension of a symmetric product. This proves the first item.

From the embedding $\IP^r\isom V\subset \IP^n$ we get the (twisted) normal bundle sequence
\[0\to \ki/\ki^2(mH)\to \Omega_{\IP^n}(mH)\restr V \to \Omega_{\IP^r}(md)\to 0.\]
The global sections of $\Omega_{\IP^r}(md)$ can be either computed via the Euler sequence or read of from Bott's formula (see e.g. \cite{bott57}).
Pulling back the Euler--sequence on $\IP^n$ with $v_d$ we get
\[ 0\to v_d^*(\Omega_{\IP^n}) (md) \to \ko_{\IP^r}((m-1)d)\tensor \Sym^1 (S_d)  \to \ko_{\IP^r}(md)\to 0.\]
The map on global sections $S_{(m-1)d}\tensor S_d \to S_{md}$ is surjective and thus
\begin{align*}
 &h^0(\IP^n, \ki/\ki^2(mH)) \\
\geq\,& h^0(\IP^r, v_d^*(\Omega_{\IP^n}) (md))-h^0(\IP^r, \Omega_{\IP^r}(md))\\
=\,& (n+1) h^0(\IP^r,\ko_{\IP^r}((m-1)d)) -h^0(\IP^r,\ko_{\IP^r}(md))-h^0(\IP^r, \Omega_{\IP^r}(md))\\
=\,&\binom{d+r}r\binom{(m-1)d+r}r- \binom{md+r}r-(md-1)\binom{r+md-1}{md},
\end{align*}
 which proves \refenum{ii}.

It remains to prove \refenum{iii}. Putting together the formulas from the first two items for $m=2$,  we {compute}
\begin{align*}
 & h^0(\IP^n, \ki/\ki^2(2H))-h^0(\IP^n, \ki(2H))\\
 \geq\,& \binom{d+r}r^2 -(2d-1)\binom{r+2d-1}{2d}-\binom{\binom{d+r}r+1}{2}\\
 = \,&\binom{d+r}r^2 -(2d-1)\binom{r+2d-1}{2d}-\frac{\binom{d+r}r(\binom{d+r}r+1)}2\\
 =\,&\frac12\left(\binom{d+r}r^2-\binom{d+r}r-2(2d-1)\binom{r+2d-1}{2d}\right).\\
\end{align*}
To obtain \refenum{iii} it therefore suffices to show that
\begin{equation}\label{eq:toshow}
\binom{d+r}r^2-\binom{d+r}r-2(2d-1)\binom{r+2d-1}{2d} > \,0.
\end{equation}
For $r=1$ {this} formula reduces to $d^2-3d+2>0$ which is certainly true for $d\geq 3$. We now proceed by induction on $r$. For convenience, note that for $r=2$ and $d=3$ the expression on the left--hand side of \eqref{eq:toshow} gives 20, so we may assume $d\geq 4$ if $r=2$.
The induction step follows from the computation below, in which we apply the induction hypothesis twice (in steps 3 and 5) and use some standard identities for binomial coefficients.
\begin{align*}
 \,&\binom{d+r+1}{r+1}^2-\binom{d+r+1}{r+1}-2(2d-1)\binom{r+2d}{2d}\\
=\,& \left(\binom{d+r}r+\binom{d+r}{r+1}\right)^2 -\binom{d+r}r-\binom{d+r}{r+1}\\
\,&\qquad-2(2d-1)\left(\binom{r+2d-1}{2d}+\binom{r+2d-1}{2d-1}\right)\\
=\,&\left(\binom{d+r}r^2-\binom{d+r}r-2(2d-1)\binom{r+2d-1}{2d}\right)\\
\,& +2\binom{d+r}r\binom{d+r}{r+1}+\binom{d+r}{r+1}^2-\binom{d+r}{r+1}-2(2d-1)\binom{r+2d-1}{2d-1}\\
>\,& 2\binom{d+r}r\binom{d+r}{r+1}+\binom{d+r}{r+1}^2-\binom{d+r}{r+1}-2(2d-1)\binom{r+2d-1}{2d-1}\\
=\,&\frac{2d}{r+1}\binom{d+r}{r}^2+\left(\frac{d}{r+1}\binom{d+r}{r}\right)^2\\
&\qquad-\frac{d}{r+1}\binom{d+r}{r}-\frac{2d}{r}2(2d-1)\binom{r+2d-1}{2d}\\
>\,&\frac{d^2+2d(r+1)}{(r+1)^2}\binom{d+r}{r}^2-\frac{d}{r+1}\binom{d+r}{r}-\frac{2d}{r}\left( \binom{d+r}r^2-\binom{d+r}r \right)\\
=\,&\left(\frac{d^2+2d(r+1)}{(r+1)^2}-\frac{2d}{r}\right)\binom{d+r}{r}^2+\left(\frac{2d}r-\frac{d}{r+1}\right)\binom{d+r}{r}\\
>\,&\frac{d}{r(r+1)^2}\left(r(d-2)-2\right)\binom{d+r}{r}^2\\
>\,&0,
\end{align*}
{where in the last step we used $d \geq 3$, and $d\geq 4$ if $r=2$.} This concludes the induction and the proof of \refenum{iii}.
\end{proof}

\subsection{Cotorsion in $\Omega_{X_{r,d}}$}
In this section we compute the cotorsion of the cones over the Veronese embeddings using their realisations as cyclic quotient singularities.

\begin{prop}\label{Vrdkt}
For all $r \geq 1, d\geq 2$ the sheaf $\Omega_{X_{r,d}}$ has cotorsion.
\end{prop}
\begin{proof}
By Corollary \ref{cotors} we need to compare the number of {generators} of $\Omega_{\IA^{r+1}}^{\mu_d}$ and the embedding dimension of $X_{r,d}$.

Recall that the ring $k[x_0, \dots, x_r]^{\mu_d}$ is generated by all monomials of total degree $d$ in the coordinates $x_0, \dots, x_r$ and thus the embedding dimension of $X_{r,d}$ is $n=\binom{r+d}d$. On the other hand, the $k[X_{r,d}]$-module $\Omega_{\IA^{r+1}}^{\mu_r}$ has a minimal system of homogeneous generators given by all products of monomials of degree $r-1$ with $dx_0, \dots, dx_r$. Subtracting the embedding dimension from the number of these generators we obtain
\begin{align*}
 (r+1)\binom{d+r-1}{d-1} - \binom{d+r}{d} &= (r+1) \frac{d}{d+r}\binom{d+r}{d}-\binom{d+r}{d}\\
                           &= \frac{d(r+1) - (d+r)}{d+r}\binom{d+r}{d}\\
                           &= \frac{r(d-1)}{d+r}\binom{d+r}{d}\\
 				&> 0,
\end{align*}
where in the last step we use $d\geq 2$. It follows that $\Omega_{X_{n,r}}$ has cotorsion.
\end{proof}

\subsection*{Acknowledgements:}
We thank Duco van Straten for an interesting discussion and  Miles Reid for many helpful comments on an earlier version of this note. {After posting of the first version of this paper on the arXiv Jonathan Wahl pointed us to \cite{wahl97}, which led to a strenghening of Proposition \ref{Vrd}.}

The first author was partly supported by the DFG--Forschergruppe 790 ''Classification of Algebraic Surfaces and Compact Complex Manifolds'' (Bay\-reuth -- Freiburg). During the preparation of the second version of the paper, he enjoyed the hospitality of the Mathematics Department at Princeton University. He gratefully acknowledges the support of the ''Eliteprogramm f\"ur Postdoktorandinnen und Postdoktoranden'' of the Baden--W\"urttemberg--Siftung. The second author gratefully acknowledges support by the DFG via the SFB/TR 45 ``Periods, moduli spaces and arithmetic of algebraic varieties`` and his Emmy--Noether project, and partial support by the Hausdorff Centre for Mathematics in Bonn.

\enlargethispage{\baselineskip}

\end{document}